\documentclass[a4paper,twoside,final]{article}
\usepackage{amssymb,amsmath}
\usepackage[latin1]{inputenc}
\usepackage[T1]{fontenc}
\usepackage{lpretex}
\usepackage[all]{xy}
\usepackage{amscd}
\usepackage{title}
\usepackage{url}
\usepackage{htmlref}

\makeatletter
\let\LB@maybe@chap\relax
\makeatother

\newcommand\per[1]{\raisebox{5pt}{$\scriptstyle p$}\!#1}
\newcommand{\wt}{\widetilde}
\newcommand\card{\symb{card}}

\begin{document}

\title[On the shape of Bruhat intervals]{On the shape of Bruhat intervals}

\author
{Anders Bj\"orner}
\address{Department of Mathematics\\
         Royal Institute of Technology\\
         SE-100~44 Stockholm, Sweden}
\email{bjorner@math.kth.se}

\author
{Torsten Ekedahl}
\address{Department of Mathematics\\
         Stockholm University\\
         SE-106~91 Stockholm, Sweden}
\email{teke@math.su.se}

\maketitle

\begin{abstract}
Let $(W,S)$ be a crystallographic Coxeter group (this includes all finite and
affine Weyl groups), and $J\subseteq S$. Let $W^J$ denote the set of minimal
coset representatives modulo the parabolic subgroup $W_J$.  For $w\in W^J$, let
$f^{w,J}_{i}$ denote the number of elements of length $i$ below $w$ in Bruhat
order on $W^J$ (notation simplified to $f^{w}_{i}$ in the case when $W^J=W$).
We show that
\begin{displaymath}
0\le i<j\le\, \ell (w)-i  \quad\Longrightarrow\quad  f^{w,J}_{i} \le  f^{w,J}_{j}.
\end{displaymath}
Furthermore, the case of equalities $f^{w}_{i} = f^{w}_{\ell({w})-i}$ for $i=1,
\ldots,k$ is characterized in terms of vanishing of coefficients in the
Kazhdan-Lusztig polynomial $P_{e,w}(q)$.

It is also shown that if $W$ is finite then the number sequence 
$f^{w}_{0},  f^{w}_{1}, \ldots, f^{w}_{\ell (w)}$ cannot grow too
rapidly. Further, in the finite case, for any given $k\ge 1$ and any $w\in W$ 
of sufficiently great length (with respect to $k$):
\begin{displaymath}
f^{w}_{\ell(w)-k} \ge  f^{w}_{\ell(w)-k+1} \ge\cdots \ge f^{w}_{\ell (w)}.
\end{displaymath}

The proofs rely for the most part on properties of the
cohomology of Kac-Moody \\ Schubert varieties.
\end{abstract}
\begin{section}{Introduction}

Let $(W,S)$ be a Coxeter system with $S$ finite. Fix a subset $J\subseteq S$ and
let $W^J$ denote the set of minimal coset representatives modulo the parabolic
subgroup $W_J=\langle J\rangle$. For $w\in W^J$, let
\begin{displaymath}
f^{w,J}_{i} := \card \{u\in W^J \mid u\le w \mbox{  and  } \ell (u)=i\}.
\end{displaymath}
In words, $f^{w,J}_{i}$ is the number of length $i$ elements in the Bruhat interval
$[e,w]^J = [e,w] \cap W^J$.  
For terminology and basic facts concerning Coxeter groups, Weyl groups, and Bruhat
order, we refer to \cite{bjoerner05::combin+coxet, humphreys90::reflec+coxet}.

For a certain class of Coxeter groups we can apply geometric methods to study
the $f^{w,J}_i$. The groups to which our methods are applicable are those 
for which the order of a product of two generators is $2$, $3$, $4$, $6$, or
$\infty$; these are usually called the \emph{crystallographic} Coxeter groups. They are 
precisely the groups that appear as Weyl groups of Kac-Moody algebras 
(cf.\ \cite[Prop.\ 3.13]{kac83::infin+lie}). 
The main purpose of this paper is to prove
the following relations and some of their ramifications.
\begin{Theorem}[A]
\label{thm1}
Let $(W,S)$ be a  crystallographic Coxeter group, $J\subseteq S$
and $w\in W^J$.  We have that
\begin{displaymath}
0\le i<j\le\, \ell (w)-i  \quad\Longrightarrow\quad  f^{w,J}_{i} \le  f^{w,J}_{j}.
\end{displaymath}
\end{Theorem}

Combining the cohomological arguments used for proving Theorem \ref{thm1} 
with a linear-algebraic argument of
Stanley \cite{stanley80::weyl+lefsc+spern} we can sharpen Theorem
\ref{thm1} to a combinatorial statement 
giving structural reasons for these inequalities.
\begin{Theorem}[B]\label{disjointchains}
Let $(W,S)$ be a crystallographic Coxeter group, $J\subseteq S$.  Fix $w\in W^J$
and $i$ such that $0\le i<\ell(w)/2$.  Then, in $[e,w]^J$ there exist
$f^{w,J}_{i}$ pairwise disjoint chains $u_i < u_{i+1} < \cdots < u_{\ell(w)-i}$
such that $\ell(u_j)=j$.
\end{Theorem}
\medskip

The inequalities of Theorem \ref{thm1} are equivalent to the following two sets
of inequalities combined:
\begin{equation}
\label{eq1} 
f^{w,J}_{i} \le f^{w,J}_{\ell (w)-i}, \quad \mbox{ for all } i < \ell (w) /2 ,
\end{equation}
and
\begin{equation}
\label{eq2}
f^{w,J}_{0} \le  f^{w,J}_{1} \le\dots \le f^{w,J}_{\lceil \ell (w) /2 \rceil}.
\end{equation}

For the rest of this section we treat only the $J=\emptyset$ case.  Then $W^J=W$, so
for simplicity we drop ``$J$'' from the notation.  In this case 
the relations (\ref{eq1}) sharpen the inequalities 
$\sum_{i\le k} f^w_i \le \sum_{i\le k} f^w_{\ell(w)-i}$, for all $0\le k < \ell(w)/2$,  due to
Brion \cite[Cor.~2]{brion00::poinc}.

The case of equality in some of the relations (\ref{eq1}) is interesting. 
Fix $w\in W$, and let $m
:=\lfloor (\ell({w})-1)/2\rfloor$.
Let $$P_{e,w}(q)=1+ \beta_0 + \beta_1 q +\cdots + \beta_m q^{m}
$$
be the Kazhdan-Lusztig polynomial of the interval $[e,w]$.  It is known
\cite[Thm. 12.2.9]{kumar02::kac+moody} that all Kazhdan-Lusztig polynomials have
non-negative coefficients if $W$ is crystallographic, and that $\beta_0=0$.

\begin{Theorem}[C]\label{equalitycase}
Suppose that $(W,S)$ is crystallographic. Let $w\in W$ and $0\le k\le m$.
Then the following conditions are equivalent:
\begin{itemize}
\item[(a)] $f^{w}_{i} = f^{w}_{\ell({w})-i}$, \, for $i=0,\dots,k$,
\item[(b)] $\beta_i =0$, \, for $i=0,\dots,k$.
\end{itemize}\noindent
Furthermore, if $k<m$ then (a) and (b) imply
\begin{itemize}
\item[(c)] $\beta_{k+1}= f^{w}_{\ell({w})-k-1} - f^{w}_{k+1}$.
\end{itemize}
\end{Theorem}

\noindent In the case $k=m$ the equivalence of (a) and (b) specializes to a 
criterion for rational smoothness of the Schubert variety $X_w$ due to
Carrell and Peterson \cite{carrell94::bruhat+coxet+deodh+schub}.


The next result shows, among other things, that for finite groups  the increasing sequence 
(\ref{eq2}) cannot grow too fast. The condition of being an $M$-sequence is
recalled in Section \ref{sec:thmD}.

\begin{Theorem}[D]\label{thmD} 
Let $(W,S)$ be a finite Weyl group and $w\in W$. Then the vectors $(f^{w}_{0},
f^{w}_{1}, \ldots, f^{w}_{\ell (w)})$ and
$(f^{w}_{0},f^{w}_{1}-f^{w}_{0},f^{w}_{2}-f^{w}_{1}, \ldots, f^{w}_{\lfloor\ell
(w)/2\rfloor}-f^{w}_{\lfloor\ell (w)/2\rfloor-1})$ are $M$-sequences.
\end{Theorem}

The increasing inequalities (\ref{eq2}) have decreasing counterparts at the upper
end of the Bruhat interval, but the information we are able to give about this
is much weaker.
\begin{Theorem}[E]
\label{thm2}
For all $k\ge1$ there exists a number $N_k$, such that for every finite Coxeter  group
$(W,S)$ and every 
$w\in W$ such that $\ell (w) \ge N_k$ we have that
\begin{displaymath}
f^{w}_{\ell(w)-k} \ge  f^{w}_{\ell(w)-k+1} \ge\dots \ge f^{w}_{\ell (w)}.
\end{displaymath}
\end{Theorem}

The paper is organized as follows. 
Sections 2 and 3 contain preliminary material on the algebraic geometry
underlying the proofs of Theorem \ref{thm1} and \ref{equalitycase}
in Section 4. The proofs of
Theorems \ref{disjointchains}, \ref{thmD} and \ref{thm2} can be found
in Sections 5, 6 and 7, respectively.
Section 8 expands on some algebraic geometry needed for
the proof of Theorem \ref{equalitycase}.

\end{section}

\begin{section}{The pure cohomology}

Let $F$ be an endomorphism of a graded $\Q_\ell$-vector space $V$ which is
finite dimensional in each degree. It will be said to be \Definition{of weight
$\le w$} (resp.\ \Definition{of pure weight $w$}), with respect to a positive
integer $q$, if the eigenvalues of $F$ on $V^i$ are algebraic numbers all of
whose conjugates have the same absolute value $q^{j/2}$ for some $j \le w+i$
(resp.\ $j=w+i$). A theorem of Deligne provides a large number of such vector
spaces in the following way.

Consider a proper variety $X_0$ over a finite field
$\F_q$ and \'etale cohomology $H^*(X,\Q_\ell)$, considered as a graded vector
space. (We shall follow the usual convention of using $0$ as a subscript to
denote objects over a finite field and drop the subscript when we extend scalars
to an algebraic closure of that field). The Frobenius map on $X_0$ induces an
endomorphism $F$ of $H^*(X,\Q_\ell)$ and Deligne's theorem \cite[5.1.14]{BBD82}
says that the action of $F$ on $H^*(X,\Q_\ell)$ is of weight $\le 0$. We shall
only be interested in the \Definition{pure part}, $H^*_p(X,\Q_\ell)$, of
$H^*(X,\Q_\ell)$, which by definition is obtained from $H^*(X,\Q_\ell)$ by
factoring out by the $F$-generalised eigenspaces of weight $<0$. Our first
result will identify $H^*_p(X,\Q_\ell)$ with the image of $H^*(X,\Q_\ell)$ in
the (middle perversity) \emph{intersection cohomology} $\cI H^*(X,\Q_\ell)$ of
$X$.

We start by recalling how the map from ordinary cohomology to intersection
cohomology is induced from a map of complexes of sheaves. Let $X$ be an
algebraic variety (we assume varieties to be reduced but not necessarily
irreducible) and \map{j}{U}{X} the inclusion of its non-singular locus. The
natural map $j_!\Q_\ell \to Rj_*\Q_\ell$ factors as $j_!\Q_\ell \to \Q_\ell
\to Rj_*\Q_\ell$. This gives maps $\per H^0(j_!\Q_\ell) \to \per H^0(\Q_\ell)
\to \per H^0(Rj_*\Q_\ell)$. Considering the distinguished triangle $\to
j_!\Q_\ell \to \Q_\ell \to i_*\Q_\ell \to$, where \map{i}{Y}{X} is the inclusion
of the singular locus, and as $i_*\Q_\ell \in \per D_c^{\le -1}(X)$, we get that
$\per H^0(j_!\Q_\ell) \to \per H^0(\Q_\ell)$ is surjective. As $j_{!*}\Q_\ell$
(by definition) is the image of $\per H^0(j_!\Q_\ell) \to \per H^0(Rj_*\Q_\ell)$
we get a surjective map $\per H^0(\Q_\ell) \to j_{!*}\Q_\ell$, and as $\Q_\ell
\in \per D^{\le 0}_c(X)$ we have a map $\Q_\ell \to \per H^0(\Q_\ell)$. So, by
composition we get a map $\Q_\ell \to j_{!*}\Q_\ell$. This map induces the
desired map $H^*(X,\Q_\ell) \to \cI H^*(X,\Q_\ell)$ from cohomology to
intersection cohomology.

Note now that, by definition, if $X_0$ is an algebraic variety over a finite
field $\F_q$, a mixed complex $\cE_0$ on $X_0$ (cf. \cite[1.1.2]{De80}) is of
weight $\le w$ precisely when the graded vector space $H^*(\cE_{\overline{s}})$
is of weight $\le w$ (with respect to $q$) for every geometric point $\overline{s}$ of $X_0$
with image a closed point of $X_0$ whose residue field has cardinality $q$.
\begin{theorem}\label{injectivity}
Let $X_0$ be a proper variety over a finite field. Then the kernel of the map
from ordinary cohomology to intersection cohomology
$H^*(X,\Q_\ell) \to \cI H^*(X,\Q_\ell)$ consists exactly of the part of
$H^*(X,\Q_\ell)$ of weight $<0$. In particular, $H^*_p(X,\Q_\ell)$ is the image
of $H^*(X,\Q_\ell)$ in $ \cI H^*(X,\Q_\ell)$.
\begin{proof}
Let \map{j}{U_0}{X_0} be the inclusion of the smooth locus.
We fit the map $\Q_\ell \to j_{!*}\Q_\ell$ into a distinguished triangle
\begin{displaymath}
\to \cF \to \Q_\ell \to j_{!*}\Q_\ell \to
\end{displaymath}
and we start by showing that $\cF$ is of weight $<0$. We know, \cite[Cor
5.4.3]{BBD82}, that $j_{!*}\Q_\ell$ is pure of weight $0$. Let now
$\overline{s}$ be a geometric point of $X_0$ with image a closed point (whose
residue field then is finite). We have an exact sequence
\begin{displaymath}
H^{i-1}((j_{!*}\Q_\ell)_{\overline{s}}) \to H^i(\cF_{\overline{s}}) \to
H^i((\Q_\ell)_{\overline{s}}) \to H^{i}((j_{!*}\Q_\ell)_{\overline{s}})
\end{displaymath}
and as $H^{i-1}((j_{!*}\Q_\ell)_{\overline{s}})$ is of weight $\le i-1$ it is
enough to show that $H^i((\Q_\ell)_{\overline{s}}) \to
H^{i}((j_{!*}\Q_\ell)_{\overline{s}})$ is injective, which is a non-trivial
condition only for $i=0$. In that case it is indeed injective as the composite
$\Q_\ell \to j_{!*}\Q_\ell \to Rj_*\Q_\ell$ induces an isomorphism on
$H^0(-)$. By Deligne's theorem \cite[5.1.14]{BBD82} we get that
$H^*(X,\cF)$ is of weight $<0$, and as the sequence
\begin{displaymath}
H^*(X,\cF) \to H^*(X,\Q_\ell) \to H^*(X,j_{!*}\Q_\ell)
\end{displaymath}
is exact, we see that the kernel of $H^*(X,\Q_\ell) \to H^*(X,j_{!*}\Q_\ell)$ is
of weight $<0$.

 On the other hand, again by Deligne's theorem and duality,
$H^*(X,j_{!*}\Q_\ell)$ is pure of weight $0$, and hence everything in
$H^*(X,\Q_\ell)$ of weight $<0$ lies in the kernel.
\end{proof}
\end{theorem}
\begin{remark}
\part Over the complex numbers this result is proved in \cite{weber::pure}
(using instead Deligne's Hodge-theoretically defined
weight filtration). For the applications of this paper
that result could also be used. In any case, the filtration by weights of
$\ell$-adic cohomology is defined for a variety over any field, commutes
appropriately under specialization of the base field, and coincides with
Deligne's Hodge-theoretic weight filtration over the complex
numbers. Thus our results are compatible with those of \cite{weber::pure}.

\part Nothing changes in the argument if one replaces $U_0$ by a smaller open dense
subset and $\Q_\ell$ by $j_*\cE$, where $\cE$ is a local system of pure weight $0$.
\end{remark}

It seems reasonable to introduce the \Definition{pure Betti numbers},
$b_i^p := \dim_{\Q_\ell}H^i_p(X,\Q_\ell)$, and we shall indeed do so. Our next
result is a weakening of some well-known numeric consequences for the Betti
numbers of a smooth and proper variety that arise as a consequence of the hard
Lefschetz theorem. The stated restriction to varieties defined over finite fields is
easily dispensed with, but we keep it to avoid too many details.
\begin{theorem}\label{inequality}
Let $X_0$ be a proper variety over a finite field, of pure dimension $n$.
We have that $b^p_i \le b^p_{i+2j}$ for all $0 \le j \le n-i$. In particular, for
$i\le n$ we have that $b^p_{n-i} \le b^p_{n+i}$.
\begin{proof}
It follows from Theorem \ref{injectivity} that the map $H^*(X,\Q_\ell) \to
\cI H^*(X,\Q_\ell)$ has $H^*_p(X,\Q_\ell)$ as its image. As it is induced by a map of
$\Q_\ell$-complexes of sheaves, this map is an $H^*(X,\Q_\ell)$-module map. In
particular, for $0 \le j \le n-i$ it commutes with multiplication by
$c_1(\cL)^j$ (the first Chern class of a line bundle),  giving a commutative diagram
\begin{displaymath}
\begin{CD}
H^{i}_p(X,\Q_\ell)@>>>\cI H^{i}(X,\Q_\ell)\\
@VV\cap c_1(\cL)^{j}V@VV\cap c_1(\cL)^{j}V\\
H^{i+2j}_p (X,\Q_\ell)@>>>\cI H^{i+2j}(X,\Q_\ell).
\end{CD}
\end{displaymath}
By what we have proven the horisontal maps are injective, and the right vertical
map is an injection by the hard Lefschetz theorem (cf. \cite[Thm
5.4.10]{BBD82}, note that there they have made a shift by $n$ of the cohomology
sheaves, hence the difference in indexing). This implies that the left vertical
map is injective, giving $$b_i^p = \dim_{\Q_\ell}H^{i}_p(X,\Q_\ell) \le
 \dim_{\Q_\ell}H^{i+2j}_p(X,\Q_\ell) = b_{i+2j}^p.$$
\end{proof}
\end{theorem}
Using these theorems for motivation and consistency we define, for any proper
variety $X$ over an algebraically closed field, the pure 
cohomology  $H^*_p(X,\Q_\ell)$ as the image
of $H^*(X,\Q_\ell)$ in $\cI H^*(X,\Q_\ell)$, and similarly for rational
coefficients when the base field is the field of complex numbers.
\end{section}

\begin{section}{The number of cells}

We define a \Definition{stratification} of a proper variety $X$
as a (necessarily finite) collection $\{V_\alpha\}_{\alpha\in I}$ of
subvarieties of $X$, called \Definition{strata}, such that $X$ is the disjoint
union of them and the closure of each stratum is a union of strata. We get a
partial order on the index set $I$ of the strata by saying that $\alpha \le
\beta$ when $V_\alpha \subseteq \overline{V_\beta}$. It then follows that if $J
\subseteq I$ is \Definition{downwards closed} (i.e., if $\alpha \le \beta$ and
$\beta \in J$ then $\alpha \in J$) then $X_J:=\cup_{\alpha \in J}V_\alpha$ is a
closed subset, as it is the union of the closed subsets
$\overline{V_\alpha}=\cup_{\beta \le \alpha}V_\beta$ for all $\alpha \in J$.

An \Definition{algebraic cell decomposition} of $X$ is a stratification for
which each stratum (which in this case will also be called a
\Definition{cell}) is isomorphic to the $n$-dimensional affine space $\A^n$
for some $n$.
\begin{theorem}\label{coh of cell decomposition}
Let $X$ be a proper variety over an algebraically closed field, having an
algebraic cell decomposition with $f_i$ cells of dimension $i$. Then:

\part  $H^{2i+1}(X,\Q_\ell)=0$ for all integers $i$. 
In particular, $b_{2i+1}=b_{2i+1}^p=0$.

\part  $H^{2i}(X,\Q_\ell)=H^{2i}_p(X,\Q_\ell)$ for all $i$, and this
space has a basis in bijection with the set of
cells of dimension $i$. In particular,
$b_{2i}=b_{2i}^p=f_i$ for all $i$.

\part Assuming that $X$ is of pure dimension $n$ we have that $f_i \le f_{j}$
for all $i \le j \le n-i$. 
\begin{proof}
For the first part, by standard specialization arguments we may assume that $X$
is defined over the algebraic closure of a finite field $\F_q$, and after
possibly extending the finite field we may assume that $X$ as well as the strata
are defined over $\F_q$.  We now prove the statements by induction over the
number of cells of the stratification, whose index set we denote by $I$. Let
$\alpha \in I$ be maximal with respect to the partial order on $I$. This means
that the cell of index $\alpha$ is an open subset $U$ of $X$ with complement
$F$, which is a proper algebraic variety with a cell decomposition with fewer
cells than $X$. By definition,  $U\cong \A^n$ for some $n$. 

We have a long exact sequence of cohomology
\begin{displaymath}
\dots\to H^i_c(U,\Q_\ell) \to H^i(X,\Q_\ell) \to H^i(F,\Q_\ell) \to\dots
\end{displaymath}
and as $X$, $U$, and $F$ are defined over $\F_q$ the Frobenius map acts on all
the vector spaces involved compatibly with all the maps. (Note that this long
exact sequence is for cohomology with compact supports, but as $X$ and $F$ are
compact, cohomology with compact supports is equal to ordinary cohomology.)
Furthermore, we have that $H^i_c(\A^n)=0$ for $i\ne 2n$ and
$H^{2n}_c(\A^n)=\Q_\ell$ is pure of weight $2n$.  For $i$ odd this shows that
$H^i(X,\Q_\ell)=0$, while for $i$ even with $i\ne 2n$ we get that
$H^i(X,\Q_\ell)=H^i(F,\Q_\ell)$. Finally, for $i=2n$ we get a short exact
sequence $\shex{H^{2n}_c(U,\Q_\ell)}{H^{2n}(X,\Q_\ell)}{H^{2n}(F,\Q_\ell)}$,
which together with the induction hypothesis and the quoted result on the
cohomology of $\A^n$ finishes the proof of part ii).

Part iii) now follows from Theorem \ref{inequality}.
\end{proof}
\end{theorem}
\end{section}

\begin{section}{Proofs of Theorems \ref{thm1} and \ref{equalitycase}}


As was mentioned in the introduction, the crystallographic groups are precisely
the ones that appear as Weyl groups of Kac-Moody algebras.  We are now going to
apply the results of previous sections to such Weyl groups. General
references for background to this material are the books by Kac
\cite{kac83::infin+lie} and Kumar \cite{kumar02::kac+moody}, for the
algebraic-geometric aspects see \cite{kumar02::kac+moody} and
\cite{slodowy86::oen+schub+kac+moody+lie}.

We start by recalling some properties of the Schubert varieties for a Kac-Moody
algebra (group). (It seems that no attempt has been made to extend the
construction of Kac-Moody groups to positive characteristic, \`a la Chevalley,
so we restrict ourselves to characteristic zero from now on.) Let $(W,S)$ be the
Weyl group of the Kac-Moody algebra (which is a Coxeter group on the generating
set $S$), pick $J \subseteq S$ and let $W_{J}$ be the subgroup of $W$ generated
by the elements of $J$. As is well-known, there is a unique element $w$ of
minimal length in any $W_{J}$-coset $\overline{w}$. The set of such elements
will be denoted $W^{J}$. For each $w \in W^{J}$ there exists (cf.
\cite[Ch. 7]{kumar02::kac+moody}, \cite[\S
2.2]{slodowy86::oen+schub+kac+moody+lie}) a complex projective variety
$\overline{X}_w$ containing locally closed subvarieties $X_{u}$ for all 
$u\in [e,w]^J$
whose closures are disjoint unions
\begin{displaymath}
\overline{X}_{u} = \biguplus_{z \le u}X_{z},
\end{displaymath}
where $z$ is assumed to be in $W^{J}$. The partial order $\le$ is the Bruhat
order, and $X_{u}$ is a subvariety of $\overline{X}_w$ isomorphic to
$\A^{\ell(u)}$ \cite[Thm.\ 2.4]{slodowy86::oen+schub+kac+moody+lie}.
\begin{remark}
The variety $\overline{X}_{w}$ depends, at least a priori, on the choice of a
dominant weight. However, we are just going to use its existence and not any
uniqueness. (It is in any case true that any two choices give varieties that are
related by algebraic maps which are homeomorphisms
\cite{slodowy86::oen+schub+kac+moody+lie}, and hence have the same cohomology.)
\end{remark}
Going further, the Kac-Moody group has a Borel subgroup $B$ which acts on each
$\overline{X}_{w}$ such that the $X_u$, $u \le w$, are the orbits. Note that $B$
is not an algebraic group but only a group scheme (i.e., not of finite
type). However, the action on any $\overline{X}_{w}$ factors through a quotient
which is an algebraic group, and we shall therefore allow ourselves to act as if
$B$ itself was an algebraic group. Using this action we get the next result.
\begin{lemma}\label{Constant on cells}
The restriction of a $B$-complex of $\overline{X}_{w}$ to some $X_u$,
for $u \le w$, is constant.
\begin{proof}
Recall (cf.\ \cite[1.8]{slodowy86::oen+schub+kac+moody+lie}) that there is a
subvariety $U_u$ of $B$ and a point $x$ on $X_u$ such that the map $g \mapsto
gx$ gives an isomorphism $U_u \to X_u$. Now, if $C$ is the $B$-complex, then by
assumption we have an isomorphism between $p_2^*C$ and $m^*C$ on $B\times
\overline{X}_{\overline{w}}$. Taking the restriction of this isomorphism to
$U_u\times\{x\}$ we obtain an isomorphism between $C$ and the constant extension
of $C_x$ to $X_u$.
\end{proof}
\end{lemma}
\begin{remark}
By Proposition \ref{G-action on middle extension} this result can be applied to
the intersection complex. In the setup of \cite{braden01::from}, rather than
assuming a $B$-action, an equisingularity condition along $X_u$ is assumed (which
should follow from the fact that $X_u$ is a $B$-orbit). This implies that the
restriction to $X_u$ is locally constant and then the fact that $X_u$ is
isomorphic to an affine space and hence contractible implies that locally
constant complexes are constant.
\end{remark}
Using the Kac-Moody Schubert varieties we can now undertake the proofs of the
main results.
\begin{proofof}{Theorem \ref{thm1}}\label{thm1proof}
This follows immediately from Proposition \ref{coh of cell decomposition}
applied to $\overline{X}_w$.
\end{proofof}
\begin{remark}
In connection with the inequalities of Theorem \ref{thm1} it might be tempting to
speculate that the $f$-vectors $(f^{w,J}_{0}, f^{w,J}_{1}, \ldots, f^{w,J}_{\ell
(w)})$ are \emph{unimodal}, meaning that they increase up to some maximum and
then decrease. However, this is false. See Stanton \cite{stanton90::uenim} for
non-unimodal examples in groups of type $A_n$ modulo maximal parabolic
subgroups.
\end{remark}
For the proof of Theorem \ref{equalitycase} we need to extend the monotonicity
theorem of Braden and MacPherson (see \cite[Cor.\ 3.7]{braden01::from}) to the
case of a general Kac-Moody Schubert variety.  

\begin{theorem}\label{monot.thm}
Let $x\le y\le z$ in a crystallographic Coxeter group. Then $P^i_{x,z} \ge
P^i_{y,z}$, where $P^i_{\cdot ,z}$ denotes the coefficients of $q^i$ of the
respective Kazhdan-Lusztig polynomials.
\begin{remark}
It is no doubt true that the Whitney stratification condition, which is one of
the standing hypotheses of \cite{braden01::from}, is indeed fulfilled also in
this case. But rather than trying to verify that, we note that the fact that the
Schubert cells are the orbits of a group action can be used more directly to
prove the necessary conditions. As the proofs of \cite{braden01::from} are also
sometimes somewhat sketchy we have therefore chosen to go through the needed
steps rather than leaving to the reader the task of checking that the proofs of
Braden and MacPherson go through. At the same time this allows us to give the
results in our context of $\ell$-adic cohomology rather than in the de
Rham-cohomology context of \cite{braden01::from}. We have deferred to Section
\ref{sec:gmaction} the part of the argument that does not directly pertain to
Kac-Moody Schubert varieties.  We will here take that material for granted.
\end{remark}
\begin{proof}
 Here, as in the more familiar case of finite and
affine Weyl groups, these Kazhdan-Lusztig
coefficients can be interpreted as the dimension of
the fibre of the cohomology of the intersection complex for $\overline{X}_z$ at
a point of $X_x$. Hence the monotonicity theorem \ref{monot.thm} 
follows from our analogue of
\cite[Thm.\ 3.6]{braden01::from}, which we now have the appropriate tools for proving.
We state it as a separate proposition.
\end{proof}
\end{theorem}
\begin{proposition}\label{Monotonicity}
For $x \le y \le z \in W^J$ we have a surjective map of $\Q_\ell$-vector spaces
$\cI\cH^*(\overline{X}_z)_x \to \cI\cH^*(\overline{X}_z)_y$, where the
$\cI\cH^*(\overline{X}_z)$ are the cohomology sheaves of the intersection
complex of $\overline{X}_z$ and $(-)_t$ denotes a fibre at any point of $X_t$.
\begin{proof}
This is \cite[Thm.\ 3.6]{braden01::from} in our context. Its proof, 
as well as proofs of the
supporting Lemma 3.1, Lemma 3.3, Proposition 3.4, and Theorem 3.5
of \cite{braden01::from}, can now be carried through:
\begin{itemize}
\item Lemma 3.1 follows from Proposition \ref{Contracting action} and
Proposition \ref{G-action on middle extension}.

\item  Lemma 3.3, Proposition 3.4, and Theorem 3.5 can be proved with the same
proofs, using \cite{De80} instead of \cite{saito89::introd+hodge} for the weight
results (note that in sheaf theory the relative cohomology $H^*(X,U,\cF)$ is
defined as the cohomology with support $H^*_Y(X,\cF)$, where $Y:=X\setminus U$).

\item In the proof of (the analogue of) Theorem 3.6 we use Lemma \ref{Constant
on cells} to conclude that the intersection complex is constant on each $X_t$
and hence the cohomology of it on $X_t$ is equal to its fibre for any point on $X_t$.

\item To get a contracting action on $X_z$ with the same fixed points as for the
torus $T \subset B$ we use Corollary \ref{Contractible torus action} applied to
$Z$ consisting of a single fixed point (note that $\overline{X_z}$ is
irreducible). To verify the needed conditions we use that $\overline{X_z}$ lies
an ind-variety which set-theoretically is equal to $G/P_J$, where $G$ is the
Kac-Moody group and $P_J$ is the parabolic subgroup corresponding to $J$. The
fixed points of $T$ on $G/P_J$ are the cosets $wP_J$, where $w$ runs over $W^J$
(where $W$ is identified with the normaliser $N$ of $T$ in $G$ divided by $T$
itself). What we want to show is that the cone spanned by the weights of the
cotangent space of $x=xP_J \in \overline{X}_z$ does not contain a line. For this
it is enough to show the same thing for the cotangent space of $x$ in
$G/P_J$. By multiplying by $x^{-1} \in N$ we reduce to $x=e$ and by (vector
space) duality to the case of the tangent space of $e \in G/P_J$. However, that
tangent space is equal to $\frak g/\frak p_J$, where $\frak g$ is the
corresponding Kac-Moody Lie algebra and $\frak p_J$ the corresponding parabolic
subalgebra. That means that the weights form a subset of the weights of $\frak
g/\frak b$, where $\frak b$ is the Lie algebra of $B$. These weights are exactly
the negative roots of the root system of the Lie algebra and they lie in the
cone generated by the negatives of the simple roots, a cone that indeed does not
contain a line.
\end{itemize}
\end{proof}
\end{proposition}

We are now ready to prove our third main theorem.
\begin{proofof}{Theorem \ref{equalitycase}}\label{thmBproof}
Assume condition (b), and let 
\begin{displaymath}
F_w(q)= \sum_{i=0}^{\ell(w)} a_i q^{i}
:=\sum_{x\le w} q^{\ell(x)} P_{x,w} (q).
\end{displaymath}
Theorem \ref{monot.thm} implies that the
$q^{i}$-coefficient of $P_{x,w}(q)$ is $0$ for all $i=1,\ldots,k$ and all $x\le
w$.  Hence, since also $\deg P_{x,w}(q) \le \lfloor
(\ell(w)-\ell(x)-1)/2\rfloor$, we get
\begin{displaymath}
\begin{array}{lll}
a_0 =1 && a_{\ell (w)}=1 \\
a_1=f^w_1 &&  a_{\ell (w)-1}=f^w_{\ell (w)-1} \\
\vdots && \vdots  \\
a_k=f^w_k && a_{\ell (w)-k}=f^w_{\ell (w)-k} \\
a_{k+1}=f^w_{k+1}+\beta_{k+1} && a_{\ell (w)-k-1}=f^w_{\ell (w)-k-1} \\
\end{array}
\end{displaymath}
Here the last row  requires that $k<m$.

Now use that $a_{i}=a_{\ell(w)-i}$ for all $i$. This is valid in all Coxeter
groups by \cite[Lemma 2.6~(v)]{kazhdan79::repres+coxet+hecke} (in this case it
is also implied by Poincar\'e duality of middle intersection cohomology of
$\overline{X}_w$). From this we conclude condition (a), as well as (c).

Finally, assume that condition (b) fails, say $d\le k$ is minimal such that
$\beta_d \neq 0$.  Applying the implication (b) $\Rightarrow$ (c) we get that
$f^{w}_{\ell({w})-d} - f^{w}_{d}=\beta_{d}\neq 0$, so also condition (a) fails.
\end{proofof}

\end{section}

\begin{section}{Proof of Theorem \ref{disjointchains}}

We proceed with the proof of Theorem \ref{disjointchains}. As mentioned in the
introduction we shall use an argument of Stanley. When adapting it, it is both
from a geometric and linear algebra point of view more natural to consider
homology than cohomology (it will be clear that formally this is not
required). Even though there is a sheaf theoretic definition, for our purposes
it is enough to define the homology $H_i(X,\Q_\l)$, for a projective variety $X$,
as the dual vector space to $H^i(X,\Q_\l)$. Then the homology
becomes a covariant instead of contravariant functor. We denote by $f_*$ the map
induced by a map \map{f}{X}{Y}. Furthermore, if $X$ is $n$-dimensional, then
there is the \emph{trace map} $H^{2n}(X,\Q_\l) \to \Q_\l$ which is surjective
and hence gives an element, the \Definition{fundamental class}, $[X] \in
H_{2n}(X)$. If $X$ is a closed subvariety of $Y$, then we get a fundamental
class $[Y]=i_*[Y]$, the \Definition{class} of $Y$, where \map{i}{X}{Y} is the
inclusion. Now, if $U \subseteq X$ is an open $n$-dimensional subset, then the
composite $H^{2n}_c(U) \to H^{2n}(X) \to \Q_\l$ is surjective and depends only
on $U$. When $U=\A^n$, the map $H^*_c(\A^n) \to \Q_\l$ is an
isomorphism. Assuming now that $X$ has a cell decomposition, combining this with
Proposition \ref{coh of cell decomposition} (or rather its proof) gives that
$H_{2i}(X)$ has the classes of the closures of the $i$-dimensional cells as
basis.
\begin{remark}
This result may seem to contradict Theorem 2.1 of \cite{stanley80::weyl+lefsc+spern}, as
together they would imply that if $X$ is $n$-dimensional, then the number of
$i$-cells is equal to the number of $(n-i)$-cells, which is not true in general
(e.g., for a Schubert variety this is true, by the Carrell-Peterson criterion,
if and only if it is rationally smooth). When $X$ is smooth the Theorem 2.1 is
true however, and that is the only case considered in
\cite{stanley80::weyl+lefsc+spern}.
\end{remark}
Now, as $H^*(X)$ is a $\Q_\l$-algebra, $H_*(X)$ becomes a module over $H^*(X)$.
Furthermore, if \map{f}{X}{Y} is a map, then $H_*(X)$ becomes a $H^*(Y)$-module through
$y\cdot x=f^*y\cdot x$ and then, purely formally, we have the
\Definition{projection formula} $y\cdot f_*x=f_*(f^*y\cdot x)$ for $y \in H^*(Y)$ and
$x \in H_*(X)$. We are now ready to prove the analogue of \cite[Lemma
2.2]{stanley80::weyl+lefsc+spern} (which as it stands is true only in the smooth
case).
\begin{lemma}\label{Kleiman}
Let $X$ be a variety with a cell decomposition and $\cL$ a line bundle on
$X$. Then for any cell $C$ the expansion $c_1(\cL)\cdot
[\overline{C}]=\sum_Dd_{C,D}[\overline{D}]$, where $D$ runs over the cells of
$X$, has the property that $d_{C,D}=0$ unless $D \subset \overline{C}$.
\begin{proof}
Let \map{i}{C}{X} be the inclusion and consider $[\overline{C}] \in
H_*(\overline{C})$. As $\overline{C}$ has a cell decomposition, the cells of
which are the $D$ for which $D \subseteq \overline{C}$, we get that
$i^*c_1(\cL)\cdot[\overline{C}]=\sum_Dd_{C,D}[\overline{D}] \in
H_*(\overline{C})$, where $D$ runs over the cells of $\overline{C}$. Applying
$i_*$ to this formula gives $c_1(\cL)\cdot[\overline{C}]=\sum_{D \subset
\overline{C}}d_{C,D}i_*[\overline{D}]$. However, $[\overline{D}] \in
H_*(\overline{C})$ is equal to $j_*[\overline{D}]$, where
\map{j}{\overline{D}}{\overline{C}} is the inclusion. Hence
$i_*[\overline{D}]=i_*j_*[\overline{D}]=(ij)_*[\overline{D}]$, but the right
hand side is by definition the class of $\overline{D}$ in $X$.
\end{proof}
\end{lemma}
We are now almost ready to adapt the proof of
Stanley. However, in \cite{stanley80::weyl+lefsc+spern} the hard Lefschetz
theorem is combined with Lemma 1.1 for the desired conclusion and in our
situation the hard Lefschetz theorem does not quite give a bijective
map. Luckily the proof of \cite[Lemma 1.1]{stanley80::weyl+lefsc+spern} 
needs only a slight modification to be
applicable to our situation (where we also, contrary to
\cite{stanley80::weyl+lefsc+spern}, do not turn the geometric poset ``upside down'').
\begin{lemma}\label{posetchains}
Let $P$ be a finite graded poset of rank $n$. Let $P_j$ denote the set of its
elements of rank $j$ and let $V_j$ be the vector space with basis
$P_j$ over some given field. For each $i\le j< n-i$ assume given a linear
transformation \map{\varphi_j}{V_{j+1}}{V_{j}} such that the following two conditions
are satisfied:
\begin{enumerate}
\item[(a)] The composite transformation 
$\varphi_{i}\circ\varphi_{i+1}\circ\cdots\circ\varphi_{n-i-2}\circ\varphi_{n-i-1}$
is surjective.

\item[(b)] If $x\in P_{j+1}$ and $\varphi_{j}(x)=\sum_{y\in P_{j}} c^j_{x,y} y$, then
$c^j_{x,y}=0$ unless $y<x$.
\end{enumerate}\noindent
Then, in $P$ there exist $\card(P_{i})$ pairwise disjoint chains
$x_i < x_{i+1} < \cdots < x_{n-i}$  such that rank$(x_j)=j$ for all $j$.
\begin{proof}
The proof of \cite[Lemma 1.1]{stanley80::weyl+lefsc+spern} goes through in this
situation with a slight modification, for the reader's convenience we repeat the
argument. Let $m:=\card(P_i)$ and put
$\varphi:=\varphi_{i}\circ\cdots\circ\varphi_{n-i-1}$. As $\varphi$ is
surjective so is \map{\Lambda^m\varphi}{\Lambda^mV_{n-i}}{\Lambda^mV_i}. By the
definition of $\varphi$ we get that
$\Lambda^m\varphi=\Lambda^m(\varphi_i)\circ\cdots\circ\Lambda^m(\varphi_{n-i-1})$. Using the
bases $P_j$ we get bases for $\Lambda^mV_j$ and hence a matrix for each
$\Lambda^m(\varphi_j)$. An entry of the product matrix of the composite
$\Lambda^m(\varphi_i)\circ\cdots\circ\Lambda^m(\varphi_{n-i-1})$ has the form
\begin{displaymath}
\sum \det\varphi_i[Q_{i+1},Q_i]\det\varphi_{i+1}[Q_{i+2},Q_{i+1}]\cdots\det\varphi_{n-i-1}[Q_{n-i},Q_{n-i-1}],
\end{displaymath}
where $Q_{j} \subseteq P_{j}$ with $\card(Q_j)=m=\card(P_i)$, $Q_{n-i}$
specifies the entry of the matrix of the composite, $\varphi_j[Q_{j+1},Q_j]$ is
the submatrix of $\varphi_j$ corresponding to the sets of basis elements $Q_j$
and $Q_{j+1}$, and the sum runs over all the choices of $Q_j$, $i\le j<n-i$. By
assumption there is a $Q_{n-i}$ such that this sum is non-zero, and hence there
is a summand that is non-zero. This gives us a set of $Q_j$ such that all the
$\det\varphi_j[Q_{j+1},Q_{j}]$ are non-zero. In particular, one term of the
expansion of this matrix must be non-zero, which gives us a bijection
\map{\sigma_j}{Q_{j+1}}{Q_{j}} such that $c^j_{x,\sigma_j(x)}\ne 0$ for all $x
\in Q_{j+1}$. By assumption (b) this implies that $\sigma_j(x)<x$ for all $x \in
Q_{j+1}$.
\end{proof}
\end{lemma}
We can now prove Theorem \ref{disjointchains}.
\begin{proofof}{Theorem \ref{disjointchains}}
We apply Lemma \ref{posetchains} to the interval $[e,w]^J$, a graded poset of
rank $\ell(w)$. We identify the vector space $V_j$ with $H_{2j}(X)$, where $X$
is the Schubert variety of $w$ and let \map{\varphi_j}{H_{2j+2}(X)}{H_{2j}(X)}
be multiplication by $c_1(\cL)$. Condition (a) of the lemma follows from the
proof of Proposition \ref{inequality}, which shows that multiplication by
$c_1(\cL)^{n-2i}$ gives an injective map from $H^{2i}(X)$ to $H^{2n-2i}(X)$ and
hence by duality it gives a surjective map from $H_{2n-2i}(X)$ to
$H_{2i}(X)$. Condition (b) follows from Lemma \ref{Kleiman}.
\end{proofof}
\end{section}

\begin{section}{Proof of Theorem \ref{thmD}}\label{sec:thmD}

We begin by recalling the definition of an $M$-sequence.
For $n,k\ge1$ there is a unique expansion
$$
n=\binom{a_k}{k}+\binom{a_{k-1}}{k-1}+\dots+\binom{a_i}{i},
$$
with $a_k>a_{k-1}>\dots>a_i\ge i\ge1$. This given, let
\begin{eqnarray*}
\partial^k(n) &:=& \binom{a_k-1}{k-1}+\binom{a_{k-1}-1}{k-2}+\dots+
\binom{a_i-1}{i-1},\\
\partial^k(0) &:=&0.
\end{eqnarray*}
\begin{Theorem}
\textrm{(Macaulay-Stanley \cite[Thm.~2.2]{stanley78::hilber})}

For an integer sequence $(1, m_1, m_2, \ldots)$ the following conditions are
equivalent {(and this defines an $M$-sequence)}:
\begin{enumerate}
\item[(1)] $\partial^k(m_k)\le m_{k-1}$, for all $k\ge 1$,
\item[(2)] some family of monomials, closed under divisibility, contains exactly 
$m_k$ monomials of degree $k$,
\item[(3)] $\dim(A_k)=m_k$ for some graded commutative algebra
$A=\oplus_{k\ge 0} A_k$ (over some field),  such that $A$ is generated by $A_1$.
\end{enumerate}
\end{Theorem}

\begin{proofof}{Theorem \ref{thmD}}
The $f$-vector $f^{w} =\{f_0,f_1, \ldots ,f_{\ell({w})}\}$ satisfies $f_k =\dim
H^{2k}(X_w)$.  Hence, to prove that $f^{w}$ is an $M$-sequence we need to show
that $H^*(X_w)$, the cohomology algebra of the Schubert variety $X_w$ (over
$\C$), is generated in degree one (or, equivalently, in $\dim=2$).  For $w=w_0$
this is classical --- it can be seen either from the description of
$H^*(X_{w_0})$ in terms of special Schubert classes, or from the isomorphism of
$H^*(X_{w_0})$ with the coinvariant algebra of $W$.

For $w\neq w_0$ we use that the inclusion $X_w \hookrightarrow X_{w_0}$
induces an injective map on homology $H_{*}(X_w) \rightarrow H_{*}(X_{w_0})$,
as is apparent from the cell decomposition. Hence, dually there is algebra surjection
$H^*(X_{w_0}) \rightarrow H^*(X_{w})$. Since $H^*(X_{w_0})$ is generated
in degree one, so is $H^*(X_{w})$.

Finally, to prove that $(f^{w}_{0},f^{w}_{1}-f^{w}_{0},f^{w}_{2}-f^{w}_{1}, \ldots, f^{w}_{\lfloor\ell
(w)/2\rfloor}-f^{w}_{\lfloor\ell (w)/2\rfloor-1})$ is an $M$-sequence we apply
the Macaulay-Stanley  theorem to the algebra $H^*(X,\Q_\l)/c_1(\Cal
L)H^*(X,\Q_\l)+H^{>\l(w)}(X,\Q_\l)$, which by Theorem \ref{inequality}
corresponds to the desired vector.
\end{proofof}
\begin{remark}
\part
The $M$-sequence property for $f$-vectors of lower intervals $[e,w]$
does not seem to extend beyond the
case of finite groups. E.g., it  fails for the affine Weyl group $\widetilde{C}_2$,
whose Poincar\'{e} series begins
$$\sum q^{\ell(w)} = 1+3q+5q^2+8q^3+\cdots$$
Let $u\in \widetilde{C}_2$ be above all elements of length $3$. Then
$f^{u} =(1,3,5,8,\ldots)$, which is not an $M$-sequence (since
$\partial^3(8)=6 \not\le5$).

An algebraic consequence of this is that $H^*(X_u)$ is not generated in degree one.

\part The
$M$-sequence property fails also for general intervals $[x,w]$ in finite
groups.  For instance, for a particular $x\in C_4$: 
$$\sum_{x\le y\le w_0} q^{\ell(y)-\ell(x)} = 1+4q+11q^2+\cdots$$
This information can be read off from Goresky's tables
\cite{goresky81::kazhd+luszt} by letting $x=zw_0$, where $z$ is element number
377 of $C_4$.

\end{remark}

We are grateful to E. Nevo for help with finding these examples,
and to B. Shapiro for help with the proof of Theorem \ref{thmD}.

\end{section}

\begin{section}{Proof of Theorem \ref{thm2}}


Let $[u,v]$ be a Bruhat interval. The elements of $[u,v]$ of length $\ell(u)+1$
are its \emph{atoms},
We let $f^{u,v}_{\ell(u)+1}$ denote the number of atoms in $[u,v]$.

It has been shown by Dyer \cite{dyer91::oen+bruhat+coxet} that (up to
isomorphism) only finitely many posets of each given length $r$ occur as
intervals in the Bruhat order on finite Coxeter groups. Therefore, the following
function is well-defined,
\begin{displaymath}
M(r) := \mathrm{max}_{[u,v]} \; \{\, f^{u,v}_{\ell(u)+1} \mid \ell(v) -\ell(u)=r  \}.
\end{displaymath}
That is, $M(r)$ denotes the maximum number of atoms of a Bruhat interval of
length $r$ occurring in any finite Coxeter group. The initial values 
$M(2)=2,\ M(3)=4,\ M(4)=8$ are known, see \cite{hultman03::bruhat+weyl}.
Let $\wt{M}(r) :=  \mathrm{max}_{t\le r} M(t)$. (Actually,
$\wt{M}(r)=M(r)$, but we don't need this.)

As usual, denote by $w_0$ the element of maximal length in any given finite group $W$.
It follows from the classification of irreducible finite Coxeter groups
that the following number-theoretic function is well-defined,
\begin{displaymath}
Q(s):= \mathrm{max}_{(W,S)} \; \{\, \ell(w_0) \,\mid\, 
\card (W) <\infty, \: \card (S)=s  \},
\end{displaymath}
where the maximum is taken over all finite Coxeter groups of rank
$s$, of which there are only finitely many. The classification shows that
$Q(s)=s^2$ for $s\ge 9$ (the maximum occurring in type B), whereas there are
irregularities occurring for $s\le 8$ due to the exceptional groups.
\begin{lemma}\label{lem:coatoms}
Let $(W,S)$ be a finite crystallographic Coxeter group, and $w\in W$.
We have that
\begin{displaymath}
\ell(w) > Q(j) \quad\Rightarrow \quad f^w_{\ell(w)-1} > j.
\end{displaymath}
\begin{proof}
Suppose that $f^w_{\ell(w)-1} \le j$. By Theorem \ref{thm1} we have that $f^w_1
\le f^w_{\ell(w)-1} $. Hence, $f^w_1 \le j $, which means that there is a set $J
\subseteq S$ of cardinality $|J| \le j$ (the set of atoms of $[e,w]$) such that
every reduced expression for $w$ uses letters only from the set $J$ (by the
subword property of Bruhat order
\cite[Thm. 5.10]{humphreys90::reflec+coxet}). In particular, $w\in W_J$, so
$\ell(w) \le \ell(w_0 (J)) \le Q(j)$, where $w_0 (J)$ denotes the element of
maximal length in the parabolic subgroup $W_J$.
\end{proof}
\end{lemma}

\begin{proofof}{Theorem \ref{thm2}}\label{thm2proof}
Assume that $W$ is crystallographic. The easy extension of the
proof to the general case (when $W$ has components of type 
$H_3$ and $H_4$) is left to the reader.

Put $$N_k := Q(\,\wt{M}(k)-1)+k,$$
and $n:= \ell(w)\ge N_k$.
For $r$ such that $1 \le r \le k$ let $L_{n-r}:= \{x\in [e,w] \mid \ell(x)=n-r\}$.
Consider the bipartite graph with vertices $L_{n-r}\, \cup L_{n-r+1}$ and edges
$E_r = \{(x,y) \in L_{n-r} \times L_{n-r+1}\mid x<y \}$.
If $(x,y)\in E_r$ then $$\deg(x) \le \wt{M}(r) \le \wt{M}(k),$$ where
$\deg(x)$ denotes the number of edges adjacent to $x$ in $E_r$.
Similarly, by Lemma \ref{lem:coatoms},
$$\deg(y) \ge \wt{M}(k).$$ 
Thus,
$$| L_{n-r}| \cdot \wt{M}(k) \ge | E_r | \ge |L_{n-r+1}| \cdot \wt{M}(k) $$
and hence
$$f^w_{\ell(w)-r} = |L_{n-r}|  \ge |L_{n-r+1}| = f^w_{\ell(w)-r+1}.$$
\end{proofof}
In closing we would like to raise the question:
Does there exist $\alpha <1$ such that
$$f^w_{\lfloor \alpha\cdot\ell(w)\rfloor} \ge \dots \ge f^w_{\ell (w)}$$
for all $w\in W$?
\end{section}


\begin{section}{Appendix: Contracting $\mul$-actions}\label{sec:gmaction}


As was mentioned in Section 4, we have deferred to here some general
algebraic-geometric material needed for the proof of
Theorem \ref{monot.thm}.

If $(X,\cF)$ and $(Y,\cG)$ are pairs consisting of an algebraic variety and
a complex of $\l$-adic sheaves on it, then a \Definition{morphism}
$(X,\cF) \to (Y,\cG)$ consists of a map \map{f}{X}{Y} together with a map of
complexes (or more precisely a map in the derived category) $\cG \to Rf_*\cF$, or
equivalently by adjunction, a map $f^*\cG \to \cF$. Such maps form a category
with composition of $\cH \to Rg_*\cG$ and $\cG \to Rf_*\cF$ over maps $X
\mapright{f} Y \mapright{g} Z$ being obtained by first applying $Rg_*$ to $\cG
\to Rf_*\cF$ and then taking the composite $\cH \to Rg_*\cG \to
Rg_*Rf_*\cF=R(gf)_*\cF$. A morphism $(X,\cF) \to (Y,\cG)$ induces a map on
cohomology $H^*(Y,\cG) \to H^*(X,\cF)$ by applying $H^*(X,-)$ to the map $\cG
\to Rf_*\cF$ and then using $H^*(Y,Rf_*\cF)=H^*(X,\cF)$. We shall say that two
maps \map{f,g}{(X,\cF)}{(Y,\cG)} are \Definition{homotopic} if there is a
connected variety $I$, a map \map{F}{(I\times X,p_2^*\cF)}{(Y,\cG)}, where
\map{p_2}{I\times X}{X} is the projection, and two points $0,1 \in I$ such that
the composites of $F$ and the maps \map{0,1}{(X,\cF)}{(I\times X,p^*_2\cF)}
equal $f$ resp.\ $g$.
\begin{proposition}
If \map{f,g}{(X,\cF)}{(Y,\cG)} are two homotopic maps, then they induce the same
map in cohomology (over an algebraically closed field $\k$).
\begin{proof}
Using the notation in the definition of homotopicity, $f^*$ and $g^*$ factor as
$0^*F^*$ and $1^*F^*$, respectively. Hence it is enough to show that
\map{0^*=1^*}{H^*(I\times X,p^*_2\cF)}{H^*(X,\cF)}. However, $H^*(I\times
X,p^*_2\cF)=H^*(X,Rp_{1*}p_2^*\cF)$ and by the projection formula
$Rp_{1*}p_2^*\cF=R\Gamma(I,\Z_\l)\Tensor\cF$. In the latter description
$0^*$ and $1^*$ are induced by tensoring $\cF$ by the maps $R\Gamma(I,\Z_\l) \to
R\Gamma(\Sp\k,\Z_\l)=\Z_\l$, the equality being true as the spectrum of an
algebraically closed field is a point from the cohomological perspective,
induced by the inclusions \map{0,1}{\Sp\k}{I}. However, these maps are
completely determined by the maps induced on zero-th cohomology $H^0(I,\Z_\l) \to
\Z_{\l}$. As $I$ is connected, the structure map $I \to \Sp\k$ induces an
isomorphism $\Z_\l=H^0(\Sp\k,\Z_\l) \to H^0(I,\Z_\l)$. As the two composites
$\Sp\k \mapright{0,1}I \to \Sp\k$ both are the identity, this implies that
\map{0,1}{\Sp\k}{I} induce the same map $H^0(I,\Z_\l) \to H^0(\Sp\k,\Z_\l)$.
\end{proof}
\end{proposition}
We now want to apply this result to the case of an equivariant complex of
sheaves for a contracting action. Thus we are dealing with an algebraic action
of a $1$-dimensional torus $\mul:=\Sp\k[t,t^{-1}]$ on an algebraic variety $X$,
i.e., a map \map{m}{\mul\times X}{X} fulfilling the required conditions for an
action. It is said to be \Definition{contracting} if it can be extended to a map
$\A^1\times X \to X$. Our next step says, roughly, that the cohomology of a
$\mul$-equivariant complex of sheaves on $X$ (assuming that the action is
contracting) is equal to the cohomology of the restriction of the complex to the
fixed point set. There are (at least) two possible definitions of what we should
mean by an equivariant complex of sheaves (the difference being whether we
demand an action on an actual complex or only on an object of the derived
category). It turns out that the complexes to which we want to apply our results
have an action of $\mul$ in the strongest possible sense, but on the other hand
to prove our results we need a very weak notion of action which is very easy to
verify. We therefore say that a complex $C$ of $\ell$-adic sheaves is a
$\mul$-complex if there exists an isomorphism \map{\phi}{p_2^*C}{m^*C} in the
derived category, where \map{p_2}{\mul\times X}{X} is the projection on the
second component and $m$ is the given action as above. Note that one has to be
very careful when using this notion; we make for instance no claim that
$\mul$-sheaves form a triangulated category.
\begin{proposition}\label{Contracting action}
Let $X$ be a variety with a contracting $\mul$-action and $C$ a $\mul$-complex
of $\l$-adic sheaves on $X$. If $Z$ is the fixed point variety of the action, then
the restriction map $H^*(X,C) \to H^*(Z,C)$ is an isomorphism.
\begin{proof}
Let \map{i}{Z}{X} and \map{j}{U}{X} be the inclusions, where $U:=X\setminus
Z$. We then have a distinguished triangle $\to j_!j^*C \to C \to i_*C \to $ which
pulls back by $m$ and $p_2$ to distinguished triangles, and the assumed
isomorphism \map{\phi}{p_2^*C}{m^*C} induces a morphism of distinguished
triangles. By the $5$-lemma we are thus reduced to the cases where either $C$ is
supported on $Z$ or zero on it. The first case is clear, so we may assume that
$C$ is zero on $Z$. We now have two maps $(X,C) \to (X,C)$, one being the identity map and
the other the identity map on $X$ and zero on the complexes. If we can show that
they are homotopic, then we conclude that $H^*(X,C)$ is zero and hence the map
to the restriction is an isomorphism. 

Now, we consider the map, also denoted
$m$, \map{m}{\A^1\times X}{X}, which extends the $\mul$-action. We want to define
a map $(\A^1\times X,p^*_2C) \to (X,C)$ over $m$, i.e., a map of complexes $m^*C
\to p^*_2C$. However, $m^*C$ and $p^*_2C$ are by assumption isomorphic over
$\mul\times X$, and as $m^*C$ is zero in its complement in $\A^1\times X$ (as
$m(\{0\}\times X) \subseteq Z$) this isomorphism extends to a map $m^*C \to
p^*_2C$. Composing with \map{1}{\Sp\k\times X}{\A^1\times X} given by the point $1$ of
$\A^1$ clearly gives the identity map, and composing with
\map{0}{\Sp\k\times X}{\A^1\times X} given by the point $0$ of $\A^1$ gives the zero
map, as the image of $\{0\}\times X$, as was just observed, is contained in
$Z$ and $C$ is zero on $Z$.
\end{proof}
\end{proposition}
In order to be able to apply these results we shall need to on the one hand give
a method for deciding when a $\mul$-action is contracting, and on the other show
that the intersection complex of a variety with a $\mul$-action is a
$\mul$-complex in our sense. We start with the first problem. 

Suppose that
$X=\Sp R$ is an affine variety with a $\mul$-action. The action immediately
translates to a grading $R=\Dsum_{i \in \Z}R_i$ so that $\lambda \in
\k^*=\mul(\k)$ acts as $\lambda\cdot r=\lambda^ir$ for $r \in R_i$. The
condition that the action be contracting is then equivalent to $R_i=0$ for
$i<0$. This implies that if $x \in X(\k)$ is a $\k$-point, then the linear action
of $\mul$ on the cotangent space $\maxid_x/\maxid_x^2$ has the property that
only non-negative weights occur. (Recall that a linear representation of $\mul$
is a direct sum of $1$-dimensional representations of the form $\lambda \mapsto
\lambda^n$, and that $n$ is the weight of that subrepresentation.) It turns out
that there is a converse to this result (where we for simplicity, and because it
is the only case we shall use, only consider the irreducible case).
\begin{proposition}
Let $X$ be an irreducible variety with a $\mul$-action, and let $Z$ be a closed
subvariety of fixed points such that for each closed point of $Z$ the action of
$\mul$ on its cotangent space has only non-negative weights. Then the action of
$\mul$ is contracting on the union of the open $\mul$-invariant affine subsets
that meet $Z$.
\begin{proof}
An extension of the $\mul$-action to an $\A^1$-action is unique if it exists, as
$\mul$ is dense in $\A^1$. Hence it is enough to show that the action extends to
any open $\mul$-invariant affine subset that meets $Z$. We may therefore assume
that $X = \Sp R$ and the $\mul$-action then corresponds to a grading $R =
\Dsum_iR_i$. We may further assume that $Z$ contains the closed point
$z$. What we want to show is that $R_i=0$ for $i<0$. 

Now, as $Z$ is irreducible
we have that $R$ embeds in the local ring $R_{\maxid_z}$, and by a theorem of
Krull we have that $\cap_n\maxid_z^n=0$ in $R_{\maxid_z}$. Hence it will be
enough to show that $R_i \subseteq \maxid_z^n$ for $i<0$ and every $n$. We do
this by induction on $n$, where the case $n=1$ is true as $z$ is a fixed point
so that $\mul$ acts trivially on $R/\maxid_z$. The case $n=2$ is then true as
the map $R_i \to \maxid_z/\maxid_z^2$ is a map from a $\mul$-representation
consisting only of negative weight representations to a representation that by
assumption contains only non-negative weight representations. Multiplication in
$R$ now gives a surjective $\mul$-equivariant map $S^n(\maxid_z/\maxid_z^2) \to
\maxid_z^n/\maxid_z^{n+1}$, which implies that $\maxid_z^n/\maxid_z^{n+1}$
consists only of non-negative weight representations. By the induction
assumption we have an induced map $R_i \to \maxid_z^n/\maxid_z^{n+1}$, which
therefore also is the zero map.
\end{proof}
\end{proposition}
\begin{remark}
The case of the action of $\mul$ on $\P^2$ given by $(x:y:z) \mapsto (x:
\lambda y:\lambda^2 z)$ is instructive. The point $(1:0:0)$ is then a fixed
point whose cotangent space has weights $1$ and $2$. The only invariant affine
open subset that contains $(1:0:0)$ is $\{x \ne 0\}$, and on it the action 
visibly extends to a $\A^1$-action; $(1:y:z) \mapsto (1:
\lambda y:\lambda^2 z)$ makes sense also for $\lambda=0$.
\end{remark}
This proposition has the following corollary which is the only consequence that
we shall actually use.
\begin{corollary}\label{Contractible torus action}
Assume that $T$ is an algebraic torus (i.e., isomorphic to $\mul^n$ for some
$n>0$) and let $\Gamma$ be its group of characters (i.e., the group of algebraic
homomorphisms $T \to \mul$). Let $X$ be a variety with an action of $T$
such that there exists an ample line bundle on $X$ with a compatible
$T$-action. Suppose $z \in X$ is a $T$-fixed point such that the cone generated
by the characters that appear in the cotangent space of $z$ does not contain a
line (or equivalently no non-zero element and its inverse). Then there is an
algebraic group homomorphism $\mul \to T$ such that $X^T=X^{\mul}$, and such that
there exists an affine $T$-invariant neighbourhood of $z$ for which the action
is contracting.
\begin{proof}
Let $\Gamma^*$ be the group of cocharacters of $T$ (i.e., algebraic group
homomorphisms $\mul \to T$) which is identified with the dual of $\Gamma$ by
pairing \map{\phi}{\mul}{T} and \map{\varphi}{T}{\mul} using the integer $n$ for
which the composite $\varphi\circ\phi$ is of the form $\lambda \mapsto
\lambda^n$. We want to apply the proposition, so what we are looking for is an
element $\phi \in \Gamma^*$ such that its fixed point set is equal to that of
$T$ and such that its weights on the cotangent space of $z$ are all
non-negative. For the first condition it is well-known that there is a finite
set of hyperplanes in $\Gamma^*$ such that if $\phi$ is not contained in any of
them, then it has the same fixed point set as $T$. For the second condition it
is clear that we are looking for a $\phi$ with non-negative values on each
character of $T$ that appears in the cotangent space of $z$. By the assumption
on the cone generated by them (and the fact that $\Gamma^*$ is the dual of
$\Gamma$) there is such a $\phi$ outside of any finite number of hyperplanes, so
that we may simultaneously fulfill both conditions. Now, by the assumption that
$X$ admits a $T$-linearised ample line bundle we may find a $T$-invariant affine
neighbourhood of $z$, and then the proposition is applicable.
\end{proof}
\end{corollary}
We are now left with establishing that the intersection complex is a
$\mul$-complex. However, we shall need the corresponding results for other
actions so we put ourselves in a somewhat more general situation. Thus, if $G$
is an algebraic group acting on a variety $X$, we say that a complex $C$ of
$\l$-adic sheaves is a \Definition{$G$-complex} if the two complexes $p_2^*C$
and $m^*C$ are isomorphic in the derived category, where \map{p_2,m}{G\times
X}{X} is the projection and $G$-action map respectively.
\begin{proposition}\label{G-action on middle extension}
Let $X$ be a variety on which the algebraic group $G$ acts. Then the
intersection complex of $X$ is a $G$-complex.
\begin{proof}
We start by noticing that $m$ is the composite of $p_2$ and the automorphism
\begin{displaymath}
\function[t]{G\times X} {G\times X}
         {(g,x)}{(g,gx)}
\end{displaymath}
so that it will be enough to show that $t^*C$ is isomorphic to $C$, where $C$ is
the pullback by $p_2$ of the intersection complex on $X$. However, as $G$ is
smooth, pullback by $p_2$ takes perverse sheaves to perverse sheaves and
pulls back $j_*\Q_\ell$ to $j'_*\Q_\ell$, resp.\ $j_!\Q_\ell$ to $j'_!\Q_\ell$.
Here $j$ is the inclusion of the smooth locus $U$ of $X$ into $X$ and
$j'$ is the inclusion of $G\times U$, which is the smooth locus of $G\times X$. 
This implies that $C$ is isomorphic to the middle extension
$j'_{!*}\Q_\ell$. Now, it is characterised in a fashion which makes it
isomorphic to its pullback by any automorphism of $G\times X$.
\end{proof}
\end{proposition}
\end{section}

\newpage

\bibliography{preamble,abbrevs,alggeom,algebra,ekedahl,combinat}
\bibliographystyle{pretex}

\end{document}